\documentclass[reqno]{amsart}
\usepackage{amsmath,amsthm,amssymb,hyperref,graphicx,mathrsfs,mathtools}  

\newtheorem{theorem}{Theorem}
\newtheorem{corollary}{Corollary}
\newtheorem{lemma}{Lemma}
\newtheorem{thmx}{Theorem}
\newtheorem{remark}{Remark}

\allowdisplaybreaks
\title[Inequalities for rational functions with prescribed poles]{Inequalities for rational functions with prescribed poles}
\author{N. A. Rather$^1$}
\author{A. Iqbal$^2$}
\author{Ishfaq Dar$^3$}
\address{$^{1,2,3}$Department of Mathematics, University of Kashmir, Srinagar-190006, India}
\email{dr.narather@gmail.com, itz.a.iqbal@gmail.com, ishfaq619@gmail.com}
\begin{document}
\maketitle
\begin{abstract}
For rational functions $R(z)=P(z)/W(z)$, where  $P$ is a polynomial of degree at the most $n$ and $W(z)=\prod_{J=1}^{n}(z-a_j)$, with $|a_j|>1,$ $j\in \{1,2,\dots,n\},$ we use simple but elegant techniques to strengthen generalizations of certain results which extend some widely known polynomial inequalities of Erd\"{o}s-Lax and Tur\'{a}n to rational functions $R$. In return these reinforced results, in the limiting case, lead to the corresponding refinements of the said polynomial inequalities. As an illustration and as an application of our results, we obtain some new improvements of the Erd\"{o}s-Lax and Tur\'{a}n type inequalities for polynomials. These improved results take into account the size of the constant term and the leading coefficient of the given polynomial. As a further factor of consideration, during the course of this paper we shall demonstrate how some recently obtained results due to S. L. Wali and W. M. Shah, [\textit{Some applications of Dubinin’s lemma to rational functions with prescribed poles,} J. Math.Anal.Appl.450(2017)769–779], could have been proved without invoking the results of Dubinin [\textit{Distortion theorems for polynomials on the circle,} Sb. Math. 191(12) (2000) 1797–1807].
\end{abstract}
\section{\textbf{Preliminaries}}
\indent Let $n$ be a positive integer, $\mathbb{C}$ be the field of complex numbers and  $\mathbb P_n$ denote the linear space of all polynomials of degree at most $n$ over $\mathbb{C}$ . For each real number $k>0$, we define the following
\begin{align*}
\mathbb{J}_{k} &:= \{z \in \mathbb{C} ~: |z|=k\}\\
\mathbb{J}_{k_-} &:= \{z \in \mathbb{C} ~: |z|<k\}\\
\mathbb{J}_{k_+} &:= \{z \in \mathbb{C} ~: |z|>k\}
\end{align*}
However, for the sake of brevity, we shall denote  $\mathbb{J}_{1}$, $\mathbb{J}_{1_-}$ and $\mathbb{J}_{1_+}$ simply by $\mathbb{J}$, $\mathbb{J}_{_-}$ and $\mathbb{J}_{_+}$ respectively.
Further for any complex valued function $f$ defined on $\mathbb{J}_{k}$, we set  $$M(f,k):=\sup_{z\in\mathbb {J}_{k}}|f(z)|,$$ \\
If $P\in \mathbb P_n$ and $P^\prime$ be its derivative, then concerning the estimate of $M(P',1)$, we have the following 
\begin{align}\label{eq1} 
M(P',1)\leq nM(P,1).
\end{align}
Inequality \eqref{eq1} is a well known result due to Bernstein \cite{b}.
It is worth mentioning that equality holds in \eqref{eq1} if and only if $P$ has all its zeros at the origin, in view of this it might be natural to investigate what happens to inequality \eqref{eq1} if we impose certain conditions on the zeros of $P$. Two earliest results belonging to this circle of ideas which stand out in terms of their impact on the subsequent work carried out since then are the following:\\  \\
\emph{ if $P\in \mathbb{P}_n$ has no zeros in the open unit disc $\mathbb{J}_{_-}$, then }
\begin{align}\label{eq2} 
M(P',1)\leq \frac{n}{2}M(P,1).
\end{align}
\noindent On the other hand,\\ \\ \emph{if $P$ is a polynomial of degree $n$, which does not vanish  in $\mathbb{J}_{_+}$, then }
\begin{align}\label{eq2a} 
M(P',1)\geq \frac{n}{2}M(P,1).
\end{align}
Inequality \eqref{eq2} was speculated by Erd\"{o}s and later verified by Lax \cite{pd}, whereas inequality \eqref{eq2a} is due to Tur\'{a}n \cite{t}.\\ 
\section{\textbf{Rational Functions}}
If a rational function is meant to be a quotient of two polynomials and if for $a_j\in\mathbb J_{_+}$, $j \in \{1,2,...,n\}$, we define
\begin{align*}
W(z) &:=\prod_{j=1}^n(z-a_j),\\
B(z) &:=\prod_{j=1}^n\big(\frac{1-\overline{a}_{j}z}{z-a_j}\big),\\ 
\mathcal{R}_n &= \mathcal{R}_{n}(a_1, a_2,..., a_n):=\bigg\{\frac{P(z)}{W(z)}:~P\in\mathbb P_n\bigg\},
\end{align*}
then $\mathcal{R}_n$ is the set of all rational functions having poles possibly at $a_1, a_2,..., a_n$ and with a finite limit at $\infty$. Moreover $B(z) \in \mathcal{R}_n$ and $|B(z)|=1$ for $z\in \mathbb J$. \\
\indent Li, Mohapatra and Rodriguez \cite{8} extended Bernstein's inequality \eqref{eq1} to rational functions in $\mathcal{R}_n$ by proving the following:
\begin{thmx}\label{thmA}
If $t_1, t_2, \dots, t_n$ be the $n$ zeros of $B(z)=\lambda$ and $s_1, s_2, \dots, s_n$ be those of $B(z)=-\lambda$; $\lambda \in \mathbb{J}$. Then for every $R\in \mathcal{R}_n$ and $z \in \mathbb{J}$
\begin{align}\label{eq3}
|R'(z)| \leq \frac{1}{2}|B'(z)|\big\{ M_1 + M_2 \big\}.
\end{align}
where $M_1 = \max\big\{ |R(t_j)|:~j=1,2,\dots,n \big\}$ and $M_2 = \max\big\{|R(s_j)|:~j=1,2,\dots,n\big\}$.	
\end{thmx}
\noindent A simple consequence of which is the following extension of inequality \eqref{eq1} to rational functions in $\mathcal{R}_n$.
\begin{remark}
For every $R\in \mathcal{R}_n$ and $z \in \mathbb{J}$
\begin{align}\label{eq4}
|R'(z)|\leq|B'(z)|M(R,1).
\end{align}
\end{remark}
\noindent For rational functions in $\mathcal{R}_n$ with restricted zeros, they also proved the following results
\begin{thmx}\label{thmB}
If $R \in \mathcal{R}_n$ has all its zeros in $\mathbb{J} \cup \mathbb{J}_{_+} $, then for $z \in \mathbb{J} $,
\begin{align}\label{eq5}
|R'(z)|\leq \frac{1}{2}|B'(z)|M(R,1).
\end{align}
\end{thmx}
\begin{thmx}\label{thmC}
If $R \in \mathcal{R}_n$ has all its zeros in $\mathbb{J} \cup \mathbb{J}_{_-} $, then for $z \in \mathbb{J} $,
\begin{align}\label{eq6}
|R'(z)| \geq \frac{1}{2}\Big\{|B'(z)|-(n-m)\Big\}|R(z)|.
\end{align}
where $m$ is the number of zeros of $R$.
\end{thmx}
\noindent Evidently, Theorem \ref{thmB} extends inequality \eqref{eq2} to rational functions in $\mathcal{R}_n$, whereas Theorem \ref{thmC} is the corresponding extension of inequality \ref{eq2a} to rational functions in $\mathcal{R}_n$.\\
 
For the class of rational functions in $\mathcal{R}_n$ having no zeros in $\mathbb{J}_{_-}$,  Aziz and Shah \cite{AzSGla} obtained the following improvement of \eqref{eq3} (Theorem \ref{thmA})
\begin{thmx}\label{thmD}
Let $R \in \mathcal{R}_n$ have all its zeros in $\mathbb{J} \cup \mathbb{J}_{_+} $. If $t_1, t_2, \dots, t_n$ be the $n$ zeros of $B(z)=\lambda$ and $s_1, s_2, \dots, s_n$ be those of $B(z)=-\lambda$; $\lambda \in \mathbb{J}$. Then for $z \in \mathbb{J}$
\begin{align}\label{eq7}
|R'(z)| \leq \frac{1}{2}|B'(z)|\big\{ M_1 ^2 + M_2 ^2 \big\}^{\frac{1}{2}}.
\end{align}
where $M_1 = \max\big\{ |R(t_j)|:~j=1,2,\dots,n \big\}$ and $M_2 = \max\big\{|R(s_j)|:~j=1,2,\dots,n\big\}$.
\end{thmx}
\noindent As a generalization of Theorem \ref{thmB}, Aziz and Zargar \cite{AzZCan} obtained the following
\begin{thmx}\label{thmE}
If $R \in \mathcal{R}_n$ has all its zeros in $\mathbb{J}_k \cup \mathbb{J}_{k_+} $, $k \geq 1$, then for $z \in \mathbb{J} $,
\begin{align}\label{eq8}
|R'(z)|\leq \frac{1}{2}\bigg\{|B'(z)| - \frac{n(k-1)}{k+1}\frac{|R(z)|^2}{\big\{M(R,1)\big\}^2}\bigg\}M(R,1).
\end{align}
\end{thmx}
\noindent Aziz and Shah \cite{AzSBal} generalized Theorem \ref{thmD} by proving the following
\begin{thmx}\label{thmF}
Let $R \in \mathcal{R}_n$ have all its zeros in $\mathbb{J}_k \cup \mathbb{J}_{k_+} $, $ k\geq 1$. If $t_1, t_2, \dots, t_n$ be the $n$ zeros of $B(z)=\lambda$ and $s_1, s_2, \dots, s_n$ be those of $B(z)=-\lambda$; $\lambda \in \mathbb{J}$. Then for $z \in \mathbb{J}$
\begin{align}\label{eq9}
|R'(z)| \leq \frac{1}{2}\bigg\{|B'(z)|^2 - \frac{2n(k-1)}{k+1}\frac{|R(z)|^2|B^{\prime}(z)|}{M_1 ^2 + M_2 ^2}  \bigg\}^{\frac{1}{2}}\big( M_1 ^2 + M_2 ^2 \big)^{\frac{1}{2}}.
\end{align}
Where $M_1 = \max\big\{ |R(t_j)|:~j=1,2,\dots,n \big\}$ and $M_2 = \max\big\{|R(s_j)|:~j=1,2,\dots,n\big\}$.
\end{thmx}

For the class of rational functions in $\mathcal{R}_n$ having no zeros in $\mathbb{J}_{k_+}$, $k \leq 1,$ Aziz and Shah \cite{AzSBal} obtained the following generalization of Theorem \ref{thmC}

\begin{thmx}\label{thmG}
If $R \in \mathcal{R}_n$ has all its zeros in $\mathbb{J}_k \cup \mathbb{J}_{k_-} $, $k \leq 1$, then for $z \in \mathbb{J} $,
\begin{align}\label{eq10}
|R'(z)| \geq \frac{1}{2}\bigg\{|B'(z)|- \frac{n(1+k)-2m}{1+k} \bigg\}|R(z)|.
\end{align}
where $m$ is the number of zeros of $R$.	
\end{thmx}
\section{\textbf{Main Results}}
In this paper, we first arrive at certain interesting improvements of Theorems \ref{thmE}, \ref{thmF} and \ref{thmG}. Then, as an application of these strengthened results, we obtain some new and elegant refinements of the well known polynomial inequalities due to Erd\"{o}s-Lax \cite{pd} and M.A. Malik \cite{newmalik}. In addition to these things, we demonstrate how some recently proved results due to S. L. Wali and W. M. Shah, \cite{SwSJma} \cite{SwSJan}, could have been obtained without appealing to the results of Osserman \cite{osser} and Dubinin \cite{Dubin}. To Sum up, we asseverate that our work besides improving certain existing estimates also furnishes relatively elementary proofs of the results obtained by  S. L. Wali and W. M. Shah in \cite{SwSJma} and \cite{SwSJan}. We begin by presenting the following refinement of Theorem \ref{thmE}

\begin{theorem}\label{tR2}
	Suppose $R \in \mathcal{R}_n$ has all its zeros in $\mathbb{J}_k \cup \mathbb{J}_{k_+}, $ $k \geq 1$; that is, $R(z)$=$\frac{P(z)}{W(z)}$ with $P(z)=\alpha_m\prod_{j=1}^{m}(z-b_j)$, $\alpha_m \neq 0,$ $m\leq n, |b_j| \geq k\geq 1,~ j=1,2,\dots,m$. Then for $z \in \mathbb{J} $,
	\begin{align}\label{tR2c}
	\big| R^\prime (z) \big| \leq	\frac{1}{2}\Bigg[ |B^\prime (z)| - \frac{n(k-1)|R(z)|^2}{(k+1)\big\{M(R,1)\big\}^2}  - \frac{2|R(z)|^2}{\big\{M(R,1)\big\}^2}\bigg\{\frac{n}{k+1}-\sum_{j=1}^{m}\frac{1}{1+|b_j|} \bigg\}\Bigg]M(R,1).
	\end{align}
\end{theorem}
\begin{remark}\label{tR2rm1}
	Equality in inequality \eqref{tR2c} occurs at $z=1$ for
	\begin{align*}
	R(z)=\bigg( \frac{z+k}{z-a} \bigg)^n, \qquad B(z)=\bigg( \frac{1-az}{z-a} \bigg)^n, \quad a>1, ~    k \geq 1.
	\end{align*}
\end{remark}
\begin{remark}\label{tR2rm2}
	Under the hypothesis in Theorem \ref{tR2}, one can easily see that $$ \frac{n}{k+1}-\sum_{j=1}^{m}\frac{1}{1+|b_j|} \geq 0 $$
	and hence Theorem \ref{tR2} improves Theorem \ref{thmE}.
\end{remark}

\noindent In the case when $k=1$, Theorem \ref{tR2} immediately leads to the following refinement of Theorem \ref{thmB}
\begin{corollary}\label{tR2co1}
	If $R \in \mathcal{R}_n$ has all its zeros in $\mathbb{J} \cup \mathbb{J}_{_+}$; that is, $R(z)$=$\frac{P(z)}{W(z)}$ with $P(z)=\alpha_m\prod_{j=1}^{m}(z-b_j)$, $\alpha_m \neq 0,$ $m\leq n, |b_j| \geq 1,~ j=1,2,\dots,m$. Then for $z \in \mathbb{J} $,
	\begin{align}\label{tR2c01c}
	\big| R^\prime (z) \big| \leq	\frac{1}{2}\Bigg[ |B^\prime (z)| -  \frac{2|R(z)|^2}{\big\{M(R,1)\big\}^2}\bigg\{\frac{n}{2}-\sum_{j=1}^{m}\frac{1}{1+|b_j|} \bigg\}\Bigg]M(R,1).
	\end{align}
\end{corollary}

\noindent Although Remark \ref{tR2rm2} shows that Theorem \ref{tR2} improves Theorem \ref{thmE}, but the bound in \eqref{tR2c} demands that all the zeros of $R(z)$ be known beforehand. Since computation of zeros is not always an easy piece of work, therefore in those situations where the zeros of $R(z)$  are unspecified, one may desire to have an improvement which instead of the zeros of $R(z)$ depends upon some coefficients of $P(z)$ in  $R(z)=\frac{P(z)}{W(z)}$. Our next result serves this purpose.

\begin{corollary}\label{tR2co2}
	Suppose $R \in \mathcal{R}_n$ has all its zeros in $\mathbb{J}_k \cup \mathbb{J}_{k_+}, $ $k \geq 1$; that is, $R(z)$=$\frac{P(z)}{W(z)}$ with $P(z)=\alpha_m\prod_{j=1}^{m}(z-b_j)=\alpha_0+\alpha_1z+\dots+\alpha_mz^m$, $\alpha_m \neq 0,$ $m\leq n, |b_j| \geq k\geq 1,~ j=1,2,\dots,m$. Then for $z \in \mathbb{J} $,
	\begin{align}\label{tR2co2c}
	\nonumber \big| R^\prime (z) \big| &\leq \frac{1}{2}\Bigg[ |B^\prime (z)| - \frac{n(k-1)|R(z)|^2}{(k+1)\big\{M(R,1)\big\}^2}  - \frac{2|R(z)|^2}{(k+1)\big\{M(R,1)\big\}^2}\bigg\{ n-m  +  \frac{|\alpha_0|-k^m|\alpha_m|}{|\alpha_0|+k^m|\alpha_m|}  \bigg\}\Bigg]M(R,1)\\
	\end{align}	
\end{corollary}
\begin{proof}
	In view of inequality \eqref{tR2c}
	\begin{align}\label{tR2co2e1}
	\nonumber \big| R^\prime (z) \big| &\leq	\frac{1}{2}\Bigg[ |B^\prime (z)| - \frac{n(k-1)|R(z)|^2}{(k+1)\big\{M(R,1)\big\}^2}  - \frac{2|R(z)|^2}{\big\{M(R,1)\big\}^2}\bigg\{\frac{n}{k+1}-\sum_{j=1}^{m}\frac{1}{1+|b_j|} \bigg\}\Bigg]M(R,1)\\
	\nonumber &=	\frac{1}{2}\Bigg[ |B^\prime (z)| - \frac{n(k-1)|R(z)|^2}{(k+1)\big\{M(R,1)\big\}^2}  - \frac{2|R(z)|^2}{\big\{M(R,1)\big\}^2}\bigg\{ \frac{n-m}{1+k}  - \frac{1}{1+k} \sum_{j=1}^{m}\frac{k-|b_j|}{1+|b_j|} \bigg\}\Bigg]M(R,1)\\
	\nonumber &\leq	\frac{1}{2}\Bigg[ |B^\prime (z)| - \frac{n(k-1)|R(z)|^2}{(k+1)\big\{M(R,1)\big\}^2}  - \frac{2|R(z)|^2}{\big\{M(R,1)\big\}^2}\bigg\{   \frac{n-m}{1+k}  - \frac{1}{1+k} \sum_{j=1}^{m}\frac{k-|b_j|}{k+|b_j|} \bigg\}\Bigg]M(R,1)\\
	&= \frac{1}{2}\Bigg[ |B^\prime (z)| - \frac{n(k-1)|R(z)|^2}{(k+1)\big\{M(R,1)\big\}^2}  - \frac{2|R(z)|^2}{\big\{M(R,1)\big\}^2}\bigg\{  \frac{n-m}{1+k}  - \frac{1}{1+k} \sum_{j=1}^{m}\frac{1-\frac{|b_j|}{k}}{1+\frac{|b_j|}{k}}  \bigg\}\Bigg]M(R,1)
	\end{align}
	Using Lemma \ref{lR2} in inequality \eqref{tR2co2e1} and noting that $\frac{|b_j|}{k} \geq 1,$ $j=1,2,\dots,m$, we get
	\begin{align*}
	\big| R^\prime (z) \big| &\leq \frac{1}{2}\Bigg[ |B^\prime (z)| - \frac{n(k-1)|R(z)|^2}{(k+1)\big\{M(R,1)\big\}^2}  - \frac{2|R(z)|^2}{(k+1)\big\{M(R,1)\big\}^2}\bigg\{  n-m  -  \frac{1-\prod_{j=1}^{m}\frac{|b_j|}{k}}{1+\prod_{j=1}^{m}\frac{|b_j|}{k}}  \bigg\}\Bigg]M(R,1)\\
	&= \frac{1}{2}\Bigg[ |B^\prime (z)| - \frac{n(k-1)|R(z)|^2}{(k+1)\big\{M(R,1)\big\}^2}  - \frac{2|R(z)|^2}{(k+1)\big\{M(R,1)\big\}^2}\bigg\{ n-m  -  \frac{k^m|\alpha_m|-|\alpha_0|}{k^m|\alpha_m|+|\alpha_0|}  \bigg\}\Bigg]M(R,1)\\
	&=  \frac{1}{2}\Bigg[ |B^\prime (z)| - \frac{n(k-1)|R(z)|^2}{(k+1)\big\{M(R,1)\big\}^2}  - \frac{2|R(z)|^2}{(k+1)\big\{M(R,1)\big\}^2}\bigg\{ n-m  +  \frac{|\alpha_0|-k^m|\alpha_m|}{|\alpha_0|+k^m|\alpha_m|}  \bigg\}\Bigg]M(R,1)
	\end{align*} 
\end{proof}
\begin{remark}\label{tR2rm3}
Evidently, under the hypothesis of Corollary \ref{tR2co2} 
	$$ n-m  +  \frac{|\alpha_0|-k^m|\alpha_m|}{|\alpha_0|+k^m|\alpha_m|} \geq 0 $$
	and therefore inequality \eqref{tR2co2c} improves Theorem \ref{thmE}.
\end{remark}

\noindent As an application of Theorem \ref{tR2}, we next present the following refinement of a result due to  M.A. Malik \cite{newmalik} concerning polynomials not vanishing in the open disk $\mathbb{J}_{k_-}, k\geq 1.$
\begin{corollary}\label{tR2co4_new}
	If $P$ be a polynomial of degree $n$, having all its zeros in $\mathbb{J}_{k}\cup \mathbb{J}_{k_+}$, $k\geq 1$; that is, $P(z)=\alpha_n\prod_{j=1}^{n}(z-b_j)=\alpha_0+\alpha_1z+\dots+\alpha_nz^n,$ $\alpha_n \neq 0,$ $|b_j| \geq k\geq 1,~ j=1,2,\dots,n$. Then for $z \in \mathbb{J} $
	\begin{align}\label{tR2co4c_new}
		|P^\prime (z)|&\leq  \frac{1}{2}\Bigg[ n -  \frac{2}{k+1} \left(  \frac{n(k-1)}{2}  +    \frac{|\alpha_0|-k^m|\alpha_m|}{|\alpha_0|+k^m|\alpha_m|} \right) \frac{|P(z)|^2}{\big\{ M(P,1) \big\}^2} \Bigg] M(P,1)
	\end{align}
\end{corollary}
\begin{proof}
	Taking $R(z)=\frac{P(z)}{W(z)}$ with $W(z)=(z-\alpha)^n,$ $\alpha >k\geq 1,$ so that $B(z)=\big( \frac{1-\alpha z}{z-\alpha} \big)^n$, from Corollary \ref{tR2co2} with $m=n$, we get 
	\begin{align}\label{tR2co4e1}
	\big| R^\prime (z) \big| &\leq \frac{1}{2}\Bigg[ |B^\prime (z)| - \frac{n(k-1)|R(z)|^2}{(k+1)\big\{M(R,1)\big\}^2}  - \frac{2|R(z)|^2}{(k+1)\big\{M(R,1)\big\}^2}\bigg\{   \frac{|\alpha_0|-k^m|\alpha_m|}{|\alpha_0|+k^m|\alpha_m|}  \bigg\}\Bigg]M(R,1)\\
	\nonumber &= \frac{1}{2}\Bigg[ |B^\prime (z)| - \frac{2}{k+1} \left(  \frac{n(k-1)}{2}  +    \frac{|\alpha_0|-k^m|\alpha_m|}{|\alpha_0|+k^m|\alpha_m|} \right) \frac{|R(z)|^2}{\big\{M(R,1)\big\}^2} \Bigg]M(R,1) 
	\end{align}
	Assuming $M(R,1)=|R(e^{i\psi_0})|=\Big|\frac{P(e^{i\psi_0})}{W(e^{i\psi_0})}\Big|$ for some $\psi_0 \in [0, 2\pi)$, inequality \eqref{tR2co4e1} can be written as 
	{\scriptsize
	\begin{align}\label{tR2co4e2}
		\bigg| \frac{P^\prime (z)}{(z-\alpha)^n} - \frac{n(z-\alpha)^{n-1}P(z)}{(z-\alpha)^{2n}}\bigg| &\leq \frac{1}{2}\Bigg[ |B^\prime (z)| - \frac{2}{k+1} \left(  \frac{n(k-1)}{2}  +    \frac{|\alpha_0|-k^m|\alpha_m|}{|\alpha_0|+k^m|\alpha_m|} \right)   \frac{|P(z)|^2}{|P(e^{i\psi_0})|^2}\bigg|\frac{e^{i\psi_0}-\alpha}{z-\alpha}\bigg|^{2n} \Bigg] \frac{|P(e^{i\psi_0})|}{|e^{i\psi_0}-\alpha|^n}
	\end{align}	}
	
	Since this is true for every $\alpha >1$, letting $\alpha \rightarrow \infty $ in inequality \eqref{tR2co4e2} and noting as  $\alpha \rightarrow \infty,~ |B^\prime(z)|\rightarrow |nz|=n$ for $z\in \mathbb{J}$, we get
	\begin{align*}
		\nonumber |P^\prime (z)| &\leq \frac{1}{2}\Bigg[ n -  \frac{2}{k+1} \left(  \frac{n(k-1)}{2}  +    \frac{|\alpha_0|-k^m|\alpha_m|}{|\alpha_0|+k^m|\alpha_m|} \right)  \frac{|P(z)|^2}{|P(e^{i\psi_0})|^2} \Bigg] |P(e^{i\psi_0})|\\
		&\leq \frac{1}{2}\Bigg[ n -  \frac{2}{k+1} \left(  \frac{n(k-1)}{2}  +    \frac{|\alpha_0|-k^m|\alpha_m|}{|\alpha_0|+k^m|\alpha_m|} \right) \frac{|P(z)|^2}{\big\{ M(P,1) \big\}^2} \Bigg] M(P,1) 
	\end{align*}	
	
\end{proof}
\noindent For $k=1$, Corollary \ref{tR2co2} yields the following refinement of Theorem \ref{thmB}

\begin{corollary}\label{tR2co3}
	If $R \in \mathcal{R}_n$ has all its zeros in $\mathbb{J} \cup \mathbb{J}_{_+}$; that is, $R(z)$=$\frac{P(z)}{W(z)}$ with $P(z)=\alpha_m\prod_{j=1}^{m}(z-b_j)=\alpha_0+\alpha_1z+\dots+\alpha_mz^m$, $\alpha_m \neq 0,$ $m\leq n, |b_j| \geq 1,~ j=1,2,\dots,m$. Then for $z \in \mathbb{J} $,
	\begin{align}\label{tR2co3c}
	\nonumber \big| R^\prime (z) \big| &\leq \frac{1}{2}\Bigg[ |B^\prime (z)|   - \frac{|R(z)|^2}{\big\{M(R,1)\big\}^2}\bigg\{ n-m  +  \frac{|\alpha_0|-|\alpha_m|}{|\alpha_0|+|\alpha_m|}  \bigg\}\Bigg]M(R,1)\\
	\end{align}	
\end{corollary}

\noindent For $k=1$ corollary \ref{tR2co4_new} reduces to the following refinement of the famous Erd\"{o}s - Lax  inequality \eqref{eq2} concerning polynomials not vanishing in the open unit disk $\mathbb{J}_{_-}$
\begin{corollary}\label{tR2co4}
If $P$ be a polynomial of degree $n$, having all its zeros in $\mathbb{J} \cup \mathbb{J}_{_+}$ ; that is, $P(z)=\alpha_n\prod_{j=1}^{n}(z-b_j)=\alpha_0+\alpha_1z+\dots+\alpha_nz^n,$ $\alpha_n \neq 0,$ $|b_j| \geq 1,~ j=1,2,\dots,n$. Then for $z \in \mathbb{J} $
\begin{align}\label{tR2co4c}
|P^\prime (z)|&\leq \frac{1}{2}\Bigg[ n -   \frac{|P(z)|^2}{\big\{ M(P,1) \big\}^2} \bigg\{\frac{|\alpha_0|-|\alpha_n|}{|\alpha_0|+|\alpha_n|}  \bigg\}\Bigg] M(P,1) 
\end{align}
\end{corollary}

As our second result, we present the following refinement of Theorem \ref{thmF}
\begin{theorem}\label{tR1}
	Suppose $R \in \mathcal{R}_n$ has all its zeros in $\mathbb{J}_k \cup \mathbb{J}_{k_+} $, $ k\geq 1$; that is, $R(z)$=$\frac{P(z)}{W(z)}$ with $P(z)=\alpha_m\prod_{j=1}^{m}(z-b_j)$, $\alpha_m \neq 0,$ $m\leq n, |b_j| \geq k\geq 1,~ j=1,2,\dots,m$.  If $t_1, t_2, \dots, t_n$ be the $n$ zeros of $B(z)=\lambda$ and $s_1, s_2, \dots, s_n$ be those of $B(z)=-\lambda$; $\lambda \in \mathbb{J}$. Then for $z \in \mathbb{J}$
\begin{align}\label{tR1c}
	\big| R^\prime (z) \big| \leq \frac{1}{2}\Bigg[ |B^\prime (z)|^2 -\frac{2n(k-1)}{k+1} \frac{|R(z)|^2|B^\prime (z)|}{M_1 ^2 + M_2 ^2 } -  \frac{4|R(z)|^2|B^\prime (z)|}{M_1 ^2 + M_2 ^2 }\bigg\{ \frac{n}{1+k} - \sum_{j=1}^{m}\frac{1}{1+|b_j|} \bigg\} \Bigg]^{\frac{1}{2}}\big( M_1 ^2 + M_2 ^2 \big)^{\frac{1}{2}}.
	\end{align}
	Where $M_1 = \max\big\{ |R(t_j)|:~j=1,2,\dots,n \big\}$ and $M_2 = \max\big\{|R(s_j)|:~j=1,2,\dots,n\big\}$.
	\end{theorem}
\begin{remark}\label{tR1rm1}
	The inequality \eqref{tR1c} is best possible and equality holds for $R(z)=B(z)+\lambda,  \lambda \in \mathbb{J}.$
\end{remark}
\begin{remark}\label{tR1rm2}
Since 	$ \frac{n}{1+k} - \sum_{j=1}^{m}\frac{1}{1+|b_j|} = \frac{n-m}{1+k} + \sum_{j=1}^{m}\bigg( \frac{1}{1+k}- \frac{1}{1+|b_j|} \bigg) = \frac{n-m}{1+k} + \sum_{j=1}^{m}\frac{|b_j|-k}{(1+k)(1+|b_j|)} \geq 0$, it follows that the inequality \eqref{tR1c} refines the inequality \eqref{eq9}; that is, Theorem \ref{tR1} improves Theorem \ref{thmF}.
\end{remark}
\noindent For $k=1$, Theorem \ref{tR1} yields the following refinement of inequality \eqref{eq7} (Theorem \ref{thmD})
\begin{corollary}\label{tR1co1}
	Let $R \in \mathcal{R}_n$ have all its zeros in $\mathbb{J} \cup \mathbb{J}_{_+} $; that is, $R(z)$=$\frac{P(z)}{W(z)}$ with $P(z)=\alpha_m\prod_{j=1}^{m}(z-b_j)$, $\alpha_m \neq 0,$ $m\leq n, |b_j| \geq 1,~ j=1,2,\dots,m$.  If $t_1, t_2, \dots, t_n$ be the $n$ zeros of $B(z)=\lambda$ and $s_1, s_2, \dots, s_n$ be those of $B(z)=-\lambda$; $\lambda \in \mathbb{J}$. Then for $z \in \mathbb{J}$
\begin{align}\label{tR1co1c}
\big| R^\prime (z) \big| \leq \frac{1}{2}\Bigg[ |B^\prime (z)|^2  -  \frac{4|R(z)|^2|B^\prime (z)|}{M_1 ^2 + M_2 ^2 }\bigg\{ \frac{n}{2} - \sum_{j=1}^{m}\frac{1}{1+|b_j|} \bigg\} \Bigg]^{\frac{1}{2}}\big( M_1 ^2 + M_2 ^2 \big)^{\frac{1}{2}}.
\end{align}
Where $M_1 = \max\big\{ |R(t_j)|:~j=1,2,\dots,n \big\}$ and $M_2 = \max\big\{|R(s_j)|:~j=1,2,\dots,n\big\}$.
\end{corollary}
\noindent Concerning the bound which takes into consideration the size of some coefficients of $P(z)$ in $R(z)=\frac{P(z)}{W(z)}$, we present the following 
\begin{corollary}\label{tR1co2}
	Suppose $R \in \mathcal{R}_n$ has all its zeros in $\mathbb{J}_k \cup \mathbb{J}_{k_+} $, $ k\geq 1$; that is, $R(z)$=$\frac{P(z)}{W(z)}$ with $P(z)=\alpha_m\prod_{j=1}^{m}(z-b_j)=\alpha_0+\alpha_1 z + \dots + \alpha_mz^m$, $\alpha_m \neq 0,$ $m\leq n, |b_j| \geq k\geq 1,~ j=1,2,\dots,m$.  If $t_1, t_2, \dots, t_n$ be the $n$ zeros of $B(z)=\lambda$ and $s_1, s_2, \dots, s_n$ be those of $B(z)=-\lambda$; $\lambda \in \mathbb{J}$. Then for $z \in \mathbb{J}$	
	\begin{align}\label{tR1co2c}
	\big| R^\prime (z) \big| \leq \frac{1}{2}\Bigg[ |B^\prime (z)|^2 -\frac{2n(k-1)}{k+1} \frac{|R(z)|^2|B^\prime (z)|}{M_1 ^2 + M_2 ^2 } -  \frac{4|R(z)|^2|B^\prime (z)|}{(M_1 ^2 + M_2 ^2)(1+k) }\bigg\{ n-m  +  \frac{|\alpha_0|-k^m|\alpha_m|}{|\alpha_0|+k^m|\alpha_m|} \bigg\} \Bigg]^{\frac{1}{2}}\big( M_1 ^2 + M_2 ^2 \big)^{\frac{1}{2}}
	\end{align}
	Where $M_1 = \max\big\{ |R(t_j)|:~j=1,2,\dots,n \big\}$ and $M_2 = \max\big\{|R(s_j)|:~j=1,2,\dots,n\big\}$.
\end{corollary}
\begin{proof}
Inequality \eqref{tR1c} can be written as 
\begin{align}\label{tR1co2e1}
		\nonumber\big| R^\prime (z) \big| &\leq \frac{1}{2}\Bigg[ |B^\prime (z)|^2 -\frac{2n(k-1)}{k+1} \frac{|R(z)|^2|B^\prime (z)|}{M_1 ^2 + M_2 ^2 } -  \frac{4|R(z)|^2|B^\prime (z)|}{M_1 ^2 + M_2 ^2 }\bigg\{ \frac{n}{1+k} - \sum_{j=1}^{m}\frac{1}{1+|b_j|} \bigg\} \Bigg]^{\frac{1}{2}}\big( M_1 ^2 + M_2 ^2 \big)^{\frac{1}{2}}\\
	\nonumber &= \frac{1}{2}\Bigg[ |B^\prime (z)|^2 -\frac{2n(k-1)}{k+1} \frac{|R(z)|^2|B^\prime (z)|}{M_1 ^2 + M_2 ^2 } -  \frac{4|R(z)|^2|B^\prime (z)|}{M_1 ^2 + M_2 ^2 }\bigg\{ \frac{n-m}{1+k}  - \frac{1}{1+k} \sum_{j=1}^{m}\frac{k-|b_j|}{1+|b_j|} \bigg\} \Bigg]^{\frac{1}{2}}\big( M_1 ^2 + M_2 ^2 \big)^{\frac{1}{2}}\\
	\nonumber &\leq \frac{1}{2}\Bigg[ |B^\prime (z)|^2 -\frac{2n(k-1)}{k+1} \frac{|R(z)|^2|B^\prime (z)|}{M_1 ^2 + M_2 ^2 } -  \frac{4|R(z)|^2|B^\prime (z)|}{M_1 ^2 + M_2 ^2 }\bigg\{ \frac{n-m}{1+k}  - \frac{1}{1+k} \sum_{j=1}^{m}\frac{k-|b_j|}{k+|b_j|} \bigg\} \Bigg]^{\frac{1}{2}}\big( M_1 ^2 + M_2 ^2 \big)^{\frac{1}{2}}\\
 &= \frac{1}{2}\Bigg[ |B^\prime (z)|^2 -\frac{2n(k-1)}{k+1} \frac{|R(z)|^2|B^\prime (z)|}{M_1 ^2 + M_2 ^2 } -  \frac{4|R(z)|^2|B^\prime (z)|}{M_1 ^2 + M_2 ^2 }\bigg\{ \frac{n-m}{1+k}  - \frac{1}{1+k} \sum_{j=1}^{m}\frac{1-\frac{|b_j|}{k}}{1+\frac{|b_j|}{k}} \bigg\} \Bigg]^{\frac{1}{2}}\big( M_1 ^2 + M_2 ^2 \big)^{\frac{1}{2}}
	\end{align}
	Since $\frac{|b_j|}{k} \geq 1,$ $j=1,2,\dots,m$, therefore in view of Lemma \ref{lR2}, inequality \eqref{tR1co2e1} gives
	\begin{align}
	\nonumber \big| R^\prime (z) \big| &\leq \frac{1}{2}\Bigg[ |B^\prime (z)|^2 -\frac{2n(k-1)}{k+1} \frac{|R(z)|^2|B^\prime (z)|}{M_1 ^2 + M_2 ^2 } -  \frac{4|R(z)|^2|B^\prime (z)|}{(M_1 ^2 + M_2 ^2)(1+k) }\bigg\{ n-m  -  \frac{1-\prod_{j=1}^{m}\frac{|b_j|}{k}}{1+\prod_{j=1}^{m}\frac{|b_j|}{k}} \bigg\} \Bigg]^{\frac{1}{2}}\big( M_1 ^2 + M_2 ^2 \big)^{\frac{1}{2}}\\
	\nonumber &= \frac{1}{2}\Bigg[ |B^\prime (z)|^2 -\frac{2n(k-1)}{k+1} \frac{|R(z)|^2|B^\prime (z)|}{M_1 ^2 + M_2 ^2 } -  \frac{4|R(z)|^2|B^\prime (z)|}{(M_1 ^2 + M_2 ^2)(1+k) }\bigg\{ n-m  -  \frac{k^m|\alpha_m|-|\alpha_0|}{k^m|\alpha_m|+|\alpha_0|} \bigg\} \Bigg]^{\frac{1}{2}}\big( M_1 ^2 + M_2 ^2 \big)^{\frac{1}{2}}\\
	\nonumber &= \frac{1}{2}\Bigg[ |B^\prime (z)|^2 -\frac{2n(k-1)}{k+1} \frac{|R(z)|^2|B^\prime (z)|}{M_1 ^2 + M_2 ^2 } -  \frac{4|R(z)|^2|B^\prime (z)|}{(M_1 ^2 + M_2 ^2)(1+k) }\bigg\{ n-m  +  \frac{|\alpha_0|-k^m|\alpha_m|}{|\alpha_0|+k^m|\alpha_m|} \bigg\} \Bigg]^{\frac{1}{2}}\big( M_1 ^2 + M_2 ^2 \big)^{\frac{1}{2}}
	\end{align}
\end{proof}
\begin{remark}\label{tR1rm3}
	Since under the hypothesis of Corollary \ref{tR1co2}, $$n-m  +  \frac{|\alpha_0|-k^m|\alpha_m|}{|\alpha_0|+k^m|\alpha_m|} \geq 0,$$ therefore it follows that the inequality \eqref{tR1co2c} constitutes a refinement of Theorem \ref{thmF}.	
\end{remark}
\begin{remark}\label{tR1rm4}
	In view of the fact  that $\alpha_0 \neq 0$ and  
	\begin{align*}
	\frac{|\alpha_0|-|\alpha_m|}{|\alpha_0|+|\alpha_m|} - \frac{\sqrt{|\alpha_0|}-\sqrt{|\alpha_m|}}{\sqrt{|\alpha_0|}} = \frac{|\alpha_0|\sqrt{|\alpha_m|}+|\alpha_m|\sqrt{|\alpha_m|}-2|\alpha_m|\sqrt{|\alpha_0|}}{(|\alpha_0|+|\alpha_m|)\sqrt{|\alpha_0|}}=\frac{\big(\sqrt{|\alpha_0|} - \sqrt{|\alpha_m|}~\big)^2 \sqrt{|\alpha_m|}}{(|\alpha_0|+|\alpha_m|)\sqrt{|\alpha_0|}} \geq 0,
	\end{align*}
 the inequality \eqref{tR1co2c}, for $k = 1$, gives 
\begin{align}\label{tR1rm4c}
	\nonumber \big| R^\prime (z) \big| &\leq \frac{1}{2}\Bigg[ |B^\prime (z)|^2  -  \frac{2|R(z)|^2|B^\prime (z)|}{M_1 ^2 + M_2 ^2 }\bigg\{ n-m  +  \frac{|\alpha_0|-|\alpha_m|}{|\alpha_0|+|\alpha_m|} \bigg\} \Bigg]^{\frac{1}{2}}\big( M_1 ^2 + M_2 ^2 \big)^{\frac{1}{2}}\\
	&= \frac{|B^\prime (z)|}{2}\Bigg[  M_1 ^2 + M_2 ^2    -  \frac{2|R(z)|^2}{|B^\prime (z)|}\bigg\{ n-m  +  \frac{|\alpha_0|-|\alpha_m|}{|\alpha_0|+|\alpha_m|} \bigg\} \Bigg]^{\frac{1}{2}}\\
	\label{tR1rm4ca} &\leq \frac{|B^\prime (z)|}{2}\Bigg[   M_1 ^2 + M_2 ^2  -  \frac{2|R(z)|^2}{|B^\prime (z)| }\bigg\{ n-m  + \frac{\sqrt{|\alpha_0|}-\sqrt{|\alpha_m|}}{\sqrt{|\alpha_0|}}  \bigg\} \Bigg]^{\frac{1}{2}}
\end{align}
Inequalities \eqref{tR1rm4c} and \eqref{tR1rm4ca} for $m=n$ respectively reduce to the results \cite[Theorem 1]{SwSJan} and \cite[Theorem 1]{SwSJma}, both due to S. L. Wali and W. M. Shah. It is however worth mentioning that S. L. Wali and W. M. Shah have obtained \cite[Theorem 1]{SwSJan}  by employing a boundary refinement of the classical Schwarz lemma due to R. Osserman \cite{osser}, while as \cite[Theorem 1]{SwSJma} has been established by invoking a lemma due to V. K. Dubinin \cite{Dubin}. Clearly, in comparison to our simple and plain sailing procedure, the approaches adopted by S. L. Wali and W. M. Shah  appear to be quite ginormous.
\end{remark}
Finally, for the class of rational functions in $\mathcal{R}_n$ having no zeros in $\mathbb{J}_{k_+}$, $k \leq 1,$ we present the following refinement of Theorem \ref{thmG}
\begin{theorem}\label{tR3}
	Suppose $R \in \mathcal{R}_n$ has all its zeros in $\mathbb{J}_k \cup \mathbb{J}_{k_-}, $ $k \leq 1$; that is, $R(z)$=$\frac{P(z)}{W(z)}$ with $P(z)=\alpha_m\prod_{j=1}^{m}(z-b_j)$, $\alpha_m \neq 0,$ $m\leq n, |b_j| \leq k\leq 1,~ j=1,2,\dots,m$. Then for $z \in \mathbb{J} $,	
	\begin{align}\label{tR3c}
	\big| R^\prime (z) \big| \geq \frac{1}{2}\bigg\{|B^\prime (z)| + \frac{2m-n(1+k)}{(1+k)} + 2\bigg( \sum_{j=1}^{m}\frac{1}{1+|b_j|}-\frac{m}{1+k}\bigg)\bigg\}\big|R(z) \big|.
	\end{align}
\end{theorem}
\begin{remark}\label{tR3rm1}
	Inequality \eqref{tR3c} is best possible and equality occurs at $z=1$ for
	\begin{align*}
	R(z) = \frac{(z+k)^m}{(z-a)^n}, \qquad
	B(z) = \bigg( \frac{1-az}{z-a} \bigg)^n, \qquad 
m\leq n,~	k \leq 1 < a.
	\end{align*}
\end{remark}
\begin{remark}\label{tR3rm2}
Since under the conditions of Theorem \ref{tR3}
$$ \sum_{j=1}^{m}\frac{1}{1+|b_j|}-\frac{m}{1+k} \geq 0, $$
therefore one can see at once that Theorem \ref{tR3} improves Theorem \ref{thmG}.
\end{remark}
\noindent For $k=1,$ Theorem \ref{tR3} reduces to the following refinement of Theorem \ref{thmC}
\begin{corollary}\label{tR3co1}
If $R \in \mathcal{R}_n$ has all its zeros in $\mathbb{J} \cup \mathbb{J}_{_-}$; that is, $R(z)$=$\frac{P(z)}{W(z)}$ with $P(z)=\alpha_m\prod_{j=1}^{m}(z-b_j)$, $\alpha_m \neq 0,$ $m\leq n, |b_j| \leq 1,~ j=1,2,\dots,m$. Then for $z \in \mathbb{J} $,	
\begin{align}\label{tR3co1c}
\big| R^\prime (z) \big| \geq \frac{1}{2}\bigg\{|B^\prime (z)| - (n-m) + 2\bigg( \sum_{j=1}^{m}\frac{1}{1+|b_j|}-\frac{m}{2}\bigg)\bigg\}\big|R(z) \big|.
\end{align}
\end{corollary}
\noindent As we pointed out earlier, in certain situations where the zeros of $R(z)$ are not known in advance, one may desire to have a bound which instead of zeros of $R(z)$, takes into consideration the size of some coefficients of $P(z)$ in $R(z)=\frac{P(z)}{W(z)}$. In this direction, we present the following 
\begin{corollary}\label{tR3co2}
Suppose $R \in \mathcal{R}_n$ has all its zeros in $\mathbb{J}_k \cup \mathbb{J}_{k_-}, $ $k \leq 1$; that is, $R(z)$=$\frac{P(z)}{W(z)}$ with $P(z)=\alpha_m\prod_{j=1}^{m}(z-b_j)=\alpha_0+\alpha_1z+\dots+\alpha_mz^m$, $\alpha_m \neq 0,$ $m\leq n, |b_j| \leq k\leq 1,~ j=1,2,\dots,m$. Then for $z \in \mathbb{J} $,	
\begin{align}\label{tR3co2c}
\big| R^\prime (z) \big| \geq\frac{1}{2}\bigg\{|B^\prime (z)| + \frac{2m-n(1+k)}{(k+1)} + \frac{2k}{k+1}\bigg( \frac{k^m|\alpha_m|-|\alpha_0|}{k^m|\alpha_m|+|\alpha_0|} \bigg)\bigg\}\big|R(z) \big| .
\end{align}	
\end{corollary}
\begin{proof}
From inequality \eqref{tR3c}, we have
\begin{align}\label{tR3co2e1}
\nonumber \big| R^\prime (z) \big| &\geq \frac{1}{2}\bigg\{|B^\prime (z)| + \frac{2m-n(1+k)}{(1+k)} + 2\bigg( \sum_{j=1}^{m}\frac{1}{1+|b_j|}-\frac{m}{1+k}\bigg)\bigg\}\big|R(z) \big|\\
\nonumber &= \frac{1}{2}\bigg\{|B^\prime (z)| + \frac{2m-n(1+k)}{(1+k)} + 2 \sum_{j=1}^{m}\bigg(\frac{1}{1+|b_j|}-\frac{1}{1+k}\bigg)\bigg\}\big|R(z) \big|\\
\nonumber &=  \frac{1}{2}\bigg\{|B^\prime (z)| + \frac{2m-n(1+k)}{(1+k)} + \frac{2k}{k+1} \sum_{j=1}^{m}\frac{k-|b_j|}{k+k|b_j|}\bigg\}\big|R(z) \big|\\
\nonumber &\geq \frac{1}{2}\bigg\{|B^\prime (z)| + \frac{2m-n(1+k)}{(1+k)} + \frac{2k}{k+1} \sum_{j=1}^{m}\frac{k-|b_j|}{k+|b_j|}\bigg\}\big|R(z) \big|\\
&= \frac{1}{2}\bigg\{|B^\prime (z)| + \frac{2m-n(1+k)}{(1+k)} + \frac{2k}{k+1} \sum_{j=1}^{m}\frac{1-|b_j|/k}{1+|b_j|/k}\bigg\}\big|R(z) \big|.
\end{align} 
 Since $ \frac{|b_j|}{k} \leq 1$, therefore by invoking Lemma \ref{lR1}, we conclude from inequality \eqref{tR3co2e1} that
 \begin{align*}
 \big| R^\prime (z) \big| &\geq \frac{1}{2}\bigg\{|B^\prime (z)| + \frac{2m-n(1+k)}{(1+k)} + \frac{2k}{k+1}\bigg( \frac{1-\prod_{j=1}^{m}\frac{|b_j|}{k}}{1+\prod_{j=1}^{m}\frac{|b_j|}{k}} \bigg)\bigg\}\big|R(z) \big|\\
 &= \frac{1}{2}\bigg\{|B^\prime (z)| + \frac{2m-n(1+k)}{(k+1)} + \frac{2k}{k+1}\bigg( \frac{k^m|\alpha_m|-|\alpha_0|}{k^m|\alpha_m|+|\alpha_0|} \bigg)\bigg\}\big|R(z) \big|
 \end{align*}
\end{proof}
\begin{remark}\label{tR3rm3}
	Since under the hypothesis of Corollary \ref{tR3co2}
	$$ k^m|\alpha_m|-|\alpha_0| \geq 0 ,$$
	therefore one can see at once that the bound in \eqref{tR3co2c} is an improvement of that in Theorem \ref{thmG}.
\end{remark}
\begin{remark}\label{tR3rm4}
	As in the case of Remark \ref{tR1rm4}, in the light of the fact
	\begin{align*}
	\frac{|\alpha_m|-|\alpha_0|}{|\alpha_m|+|\alpha_0|} - \frac{\sqrt{|\alpha_m|}-\sqrt{|\alpha_0|}}{\sqrt{|\alpha_m|}} = \frac{|\alpha_m|\sqrt{|\alpha_0|}+|\alpha_0|\sqrt{|\alpha_0|}-2|\alpha_0|\sqrt{|\alpha_m|}}{(|\alpha_m|+|\alpha_0|)\sqrt{|\alpha_m|}}=\frac{\big(\sqrt{|\alpha_m|} - \sqrt{|\alpha_0|}~\big)^2 \sqrt{|\alpha_0|}}{(|\alpha_m|+|\alpha_0|)\sqrt{|\alpha_m|}} \geq 0,
	\end{align*}
the inequality \eqref{tR3co2c} for $k=1$, gives 
 \begin{align}\label{tR3rm4e1}
 \big| R^\prime (z) \big| &\geq\frac{1}{2}\bigg\{|B^\prime (z)| -(n-m)  +  \frac{|\alpha_m|-|\alpha_0|}{|\alpha_m|+|\alpha_0|} \bigg\}\big|R(z) \big|\\
 \label{tR3rm4e1bro} &\geq \frac{1}{2}\bigg\{|B^\prime (z)| -(n-m)  +  \frac{\sqrt{|\alpha_m|}-\sqrt{|\alpha_0|}}{\sqrt{|\alpha_m|}} \bigg\}\big|R(z) \big|.
 \end{align}
  Both \eqref{tR3rm4e1} and \eqref{tR3rm4e1bro} are due to S. L. Wali and W. M. Shah, (see \cite[Theorem 2]{SwSJan} and \cite[Theorem 2]{SwSJma}). However, we once again emphasize that while obtaining these results, S. L. Wali and W. M. Shah  have employed certain results due to R. Osserman \cite{osser} and V. N. Dubinin \cite{Dubin}, whereas we have arrived at the same inequalities by using relatively simple and direct arguments with no involvement of the results of Osserman\cite{osser} or Dubinin\cite{Dubin}. 
 \end{remark}
\noindent Finally, as an application of Theorem \ref{tR3}, we present the following
\begin{corollary}\label{tR3co3}
	Suppose $P$ is a polynomial of degree $n$,  having all its zeros in $\mathbb{J}_k \cup \mathbb{J}_{k_-}, $ $k \leq 1$; that is, $P(z)=\alpha_n\prod_{j=1}^{n}(z-b_j)=\alpha_0+\alpha_1z+\dots+\alpha_nz^n$, $\alpha_n \neq 0,$  $|b_j| \leq k\leq 1,~ j=1,2,\dots,n$. Then for $z \in \mathbb{J} $,	
	\begin{align}\label{tR3co3c}
|P^\prime (z)| &\geq \frac{n}{k+1}\bigg\{1  + \frac{k}{n}\bigg( \frac{k^n|\alpha_n|-|\alpha_0|}{k^n|\alpha_n|+|\alpha_0|} \bigg)\bigg\} |P(z)| 
	\end{align}	
\end{corollary}
\begin{proof}
Taking $R(z)=\frac{P(z)}{W(z)}$ with $W(z)=(z-\alpha)^n,$ $\alpha >1,$ so that $B(z)=\big( \frac{1-\alpha z}{z-\alpha} \big)^n$, from Corollary \ref{tR3co2} with $m=n$, we get 
\begin{align}\label{tR3co3e1}
	\big| R^\prime (z) \big| &\geq\frac{1}{2}\bigg\{|B^\prime (z)| + \frac{n(1-k)}{(k+1)} + \frac{2k}{k+1}\bigg( \frac{k^n|\alpha_n|-|\alpha_0|}{k^n|\alpha_n|+|\alpha_0|} \bigg)\bigg\}\big|R(z) \big| 
\end{align}
Using the fact that $R(z)=P(z)/W(z)$, inequality \eqref{tR3co3e1} can be written as 
\begin{align}\label{tR3co3e2}
\bigg| \frac{P^\prime (z)}{(z-\alpha)^n} - \frac{n(z-\alpha)^{n-1}P(z)}{(z-\alpha)^{2n}}\bigg| &\geq \frac{1}{2}\bigg\{|B^\prime (z)| + \frac{n(1-k)}{(k+1)} + \frac{2k}{k+1}\bigg( \frac{k^n|\alpha_n|-|\alpha_0|}{k^n|\alpha_n|+|\alpha_0|} \bigg)\bigg\} \frac{|P(z)|}{|z-\alpha|^n}
\end{align}	

Since this is true for every $\alpha >1$, letting $\alpha \rightarrow \infty $ in inequality \eqref{tR3co3e2} and noting as  $\alpha \rightarrow \infty,~ |B^\prime(z)|\rightarrow |nz|=n$ for $z\in \mathbb{J}$, we get
\begin{align*}
\nonumber |P^\prime (z)| &\geq \frac{1}{2}\bigg\{n + \frac{n(1-k)}{(k+1)} + \frac{2k}{k+1}\bigg( \frac{k^n|\alpha_n|-|\alpha_0|}{k^n|\alpha_n|+|\alpha_0|} \bigg)\bigg\} |P(z)| \\
&= \frac{n}{k+1}\bigg\{1  + \frac{k}{n}\bigg( \frac{k^n|\alpha_n|-|\alpha_0|}{k^n|\alpha_n|+|\alpha_0|} \bigg)\bigg\} |P(z)| 
\end{align*}
	\end{proof}
\begin{remark}
	The Inequality \eqref{tR3co3c} is new and refines a well known polynomial inequality due to Malik \cite{newmalik}
\end{remark}
\begin{remark}
	Taking $k=1$ in inequality \eqref{tR3co3c}, we get the following refinement of the Tur\'{a}n's \cite{t} inequality \eqref{eq2a} 
	\begin{align}\label{tR3co3cdubin}
	|P^\prime (z)| &\geq \frac{n}{2}\bigg\{1  + \frac{1}{n}\bigg( \frac{|\alpha_n|-|\alpha_0|}{|\alpha_n|+|\alpha_0|} \bigg)\bigg\} |P(z)| 
	\end{align}
	Inequality \eqref{tR3co3cdubin} is due to Dubinin \cite{Dubin2007}, who has obtained it by using the boundary Schwarz lemma of R. Osserman \cite{osser}, whereas we derived this from a more general result, namely Corollary \ref{tR3co3}, which we established by using simple and elegant arguments.
	\end{remark}

\section{\textbf{Lemmas}}
\noindent For the proof of our theorems, we require following lemmas. We begin by stating the following lemma: 
\begin{lemma}\label{lR1}
 If $ \langle  x_j  \rangle _{j=1}^{\infty}$ be a sequence of real numbers such that $0 \leq x_j \leq 1,~j \in \mathbb{N}$ then 
\begin{align}\label{le1}
\sum\limits_{j=1}^{n}\frac{1-x_j}{1+x_j} \geq \frac{1-\prod\limits_{j=1}^{n}x_j}{1+\prod\limits_{j=1}^{n}x_j}, \quad \forall ~ n \in \mathbb{N}
\end{align}
\end{lemma}
\begin{proof}
To prove this result we use induction on $n$. The result is trivially true for $n=1$.\\ For $n=2$
\begin{align*}
\frac{1-x_1}{1+x_1}+\frac{1-x_2}{1+x_2} \geq \frac{1-x_1x_2}{1+x_1x_2}
\end{align*}
if
\begin{align*}
\frac{2(1-x_1x_2)}{1+x_1+x_2+x_1x_2} \geq \frac{1-x_1x_2}{1+x_1x_2},
\end{align*}
that is, if $$(1-x_1)(1-x_2)\geq 0,$$ which is true, since  $x_1,x_2 \leq 1$.
Thus the result also holds for $n=2$. Assume the result is true for $n=r \in \mathbb{N}$.
Now since  $ \prod\limits_{j=1}^{r}x_j \leq 1,$ we have
\begin{align*}
\begin{split}
\sum\limits_{j=1}^{r+1}\frac{1-x_j}{1+x_j}{}&=\sum\limits_{j=1}^{r}\frac{1-x_j}{1+x_j}+\frac{1-x_{r+1}}{1+x_{r+1}}\\
& \geq \frac{1-\prod\limits_{j=1}^{r}x_j}{1+\prod\limits_{j=1}^{r}x_j}+\frac{1-x_{r+1}}{1+x_{r+1}}  \quad \textit{(by induction hypothesis)}\\
&\geq \frac{1-\prod\limits_{j=1}^{r+1}x_j}{1+\prod\limits_{j=1}^{r+1}x_j}. \quad \textit{(by the case $n=2$)}
\end{split}
\end{align*}
 This shows that the result holds for $n=r+1$ as well. Therefore by principle of mathematical induction, it follows that the result holds for all $n \in \mathbb{N}$. That completes the proof of Lemma \ref{lR1}.
\end{proof}
\noindent A direct consequence of the above lemma is the following
\begin{lemma}\label{lR2}
	If $ \langle  x_j  \rangle _{j=1}^{\infty}$ be a sequence of real numbers such that $x_j \geq 1,~j \in \mathbb{N}$ then 
	\begin{align}\label{le2}
	\sum\limits_{j=1}^{n}\frac{1-x_j}{1+x_j} \leq \frac{1-\prod\limits_{j=1}^{n}x_j}{1+\prod\limits_{j=1}^{n}x_j}, \quad \forall ~ n \in \mathbb{N}
	\end{align}
\end{lemma}
\noindent The next two lemmas which we need are due to Li, Mohapatra and Rodrigues \cite{8}
\begin{lemma}\label{lR3}
If $\lambda\in \mathbb{J}$, then the equation $B(z)=\lambda$ has exactly n simple roots (say) $t_1,t_2,\dots,t_n$, which all lie on the unit circle $\mathbb{J}$. Further if $R \in \mathcal{R}_n$ and $z\in\mathbb{J}$, then
\begin{align}\label{le3}
B^\prime (z)R(z)- R^\prime (z)\big\{ B(z)-\lambda \big\} = \frac{B(z)}{z}\sum_{j=1}^{n}C_j R(t_j)\bigg| \frac{B(z)-\lambda}{z-t_j} \bigg|^2,
\end{align}
where $C_j=C_j(\lambda)$ is defined by 
\begin{align}\label{le4}
\frac{1}{C_j} = \sum_{\nu=1}^{n}\frac{|a_\nu|^2-1}{|t_j-a_\nu|^2} \quad \text{for $j=1,2,\dots,n.$}
\end{align}
Furthermore for $z\in\mathbb{J}$, we have
\begin{align}\label{le5}
z\frac{B^\prime (z)}{B(z)} = \sum_{j=1}^{n}C_j \bigg | \frac{B(z)-\lambda}{z-t_j}  \bigg|^2
\end{align}
and
\begin{align}\label{le6}
|B^\prime (z)| = z\frac{B^\prime (z)}{B(z)} = \sum_{j=1}^{n}\frac{|a_j|^2-1}{|z-a_j|^2}
\end{align}
\end{lemma} 
\begin{lemma}
	If $R\in \mathcal{R}_n$ and $z\in\mathbb{J}$, then
	\begin{align}\label{le7a}
		|R^\prime (z)| + |\big(R^* (z)\big)^\prime| &\leq |B^\prime (z)|M(R,1)
	\end{align}
\end{lemma}
\noindent Next we need the following lemma due to Aziz and Shah \cite{AzSGla}
\begin{lemma}\label{lR4}
	Assume $t_1,t_2,\dots,t_n$ are the roots of $B(z)=\lambda$ and $s_1,s_2,\dots,s_n$ are those of $B(z)=-\lambda$; $\lambda \in \mathbb{J}$. If $R \in \mathcal{R}_n$ and $z \in \mathbb{J}$, then
	\begin{align}\label{le7}
	|R^\prime (z)|^2 + |\big(R^* (z)\big)^\prime|^2 &\leq \frac{1}{2}|B^\prime (z)|^2\big( M_1 ^2 + M_2 ^2 \big)
	\end{align} 
	where $R^* (z)= B(z) \overline{R(\frac{1}{\overline{z}})}$,  $M_1 = \max\big\{ |R(t_j)|:~j=1,2,\dots,n \big\}$ and $M_2 = \max\big\{|R(s_j)|:~j=1,2,\dots,n\big\}$. 
\end{lemma}
\noindent Lastly, we also require the following lemma due to Aziz and Zargar \cite{AzZCan}
\begin{lemma}\label{lR5}
	If $z \in \mathbb{J}$, then
	\begin{align}\label{le8}
\Re	\bigg(\frac{zW^\prime (z)}{W(z)}\bigg) = \frac{n-|B^\prime (z)|}{2}
	\end{align}
\end{lemma}
\section{\textbf{Proof of Theorems}}

\begin{proof}[\textbf{Proof of Theorem \ref{tR2}}]
	For $z\in\mathbb{J}$ and $R(z)$=$\frac{P(z)}{W(z)}$ with $P(z)=\alpha_m\prod_{j=1}^{m}(z-b_j)$, $\alpha_m \neq 0,$ $m\leq n, |b_j| \geq k\geq 1,~ j=1,2,\dots,m$, Lemma \ref{lR5} and straight forward computations show that
	\begin{align}\label{tR1e1}
	\nonumber \Re\bigg( \frac{zR^\prime (z)}{R(z)}\bigg) &= \Re\bigg( \frac{zP^\prime (z)}{P(z)} - \frac{zW^\prime (z)}{W(z)}  \bigg)\\
	\nonumber &= \Re\bigg( \frac{zP^\prime (z)}{P(z)}  \bigg) - \Re\bigg( \frac{zW^\prime (z)}{W(z)}  \bigg)\\
	\nonumber &= \sum_{j=1}^{m}\Re\bigg( \frac{z}{z-b_j} \bigg) - \frac{n-|B^\prime (z)|}{2}\\
	& \leq \sum_{j=1}^{m}\frac{1}{1+|b_j|} - \frac{n-|B^\prime (z)|}{2}
	\end{align}
	Also for $z \in \mathbb{J}$, one can easily see that
	\begin{align*}
	\nonumber |\big(R^* (z)\big)^\prime|&=|zB^\prime (z)\overline{R(z)}-B(z)\overline{zR^\prime (z)}|\\
	\nonumber &=\bigg|\frac{zB^\prime (z)}{B(z)}\overline{R(z)}-\overline{zR^\prime (z)}\bigg| \qquad \text{$\because$ $|B(z)|=1$ for $|z|=1$}.
	\end{align*}
	This with the help of \eqref{le6} gives 
	\begin{align*}
	\nonumber \big|\big(R^* (z)\big)^\prime\big| &= \big||B^\prime (z)|\overline{R(z)}-\overline{zR^\prime (z)}\big|\\
	&=\big||B^\prime (z)|R(z)-zR^\prime (z)\big|.
	\end{align*}
	That is 
	\begin{align}\label{tR1e2}
	\nonumber \bigg| \frac{\big(R^* (z)\big)^\prime}{R(z)} \bigg|^2 
	&=\bigg||B^\prime (z)|- \frac{zR^\prime (z)}{R(z)}\bigg|^2\\
	&= |B^\prime (z)|^2 + \bigg| \frac{zR^\prime (z)}{R(z)} \bigg|^2 - 2|B^\prime (z)| \Re\bigg( \frac{zR^\prime (z)}{R(z)} \bigg).
	\end{align}
	Using inequality \eqref{tR1e1} in  \eqref{tR1e2}, we get
	\begin{align*}
	\bigg| \frac{\big(R^* (z)\big)^\prime}{R(z)} \bigg|^2 
	&\geq |B^\prime (z)|^2 + \bigg| \frac{zR^\prime (z)}{R(z)} \bigg|^2 - 2|B^\prime (z)| \sum_{j=1}^{m}\frac{1}{1+|b_j|} + |B^\prime (z)|\big(n-|B^\prime (z)|\big)\\
	&= \bigg| \frac{zR^\prime (z)}{R(z)} \bigg|^2 + n|B^\prime (z)| -2|B^\prime (z)| \sum_{j=1}^{m}\frac{1}{1+|b_j|}.
	\end{align*}
	Which for $z \in \mathbb{J}$, leads to
	\begin{align}\label{tR1e3}
	\big| \big(R^* (z)\big)^\prime \big|^2  \geq \big| R^\prime (z) \big|^2 + n|B^\prime (z)||R(z)|^2 -2|B^\prime (z)| |R(z)|^2 \sum_{j=1}^{m}\frac{1}{1+|b_j|}.
	\end{align}
	In view of the fact that 
	$$  \big| R^\prime (z) \big|^2 + n|B^\prime (z)||R(z)|^2 -2|B^\prime (z)| |R(z)|^2 \sum_{j=1}^{m}\frac{1}{1+|b_j|} \geq 0 \quad \text{for $|b_j| \geq 1, j=1,2,\dots,m,$ with $m \leq n$ }, $$
	inequality \eqref{tR1e3} together with inequality \eqref{le7a} gives
	\begin{align*}
	\big| R^\prime (z) \big| + \Bigg[\big| R^\prime (z) \big|^2 + n|B^\prime (z)||R(z)|^2 -2|B^\prime (z)| |R(z)|^2 \sum_{j=1}^{m}\frac{1}{1+|b_j|} \Bigg]^{\frac{1}{2}} &\leq \big|R^\prime (z)\big| + \big|\big(R^* (z)\big)^\prime\big|\\
	&\leq |B^\prime (z)|M(R,1).
	\end{align*}
	Which after straightforward simplifications  gives
	\begin{align*}
	\big| R^\prime (z) \big| &\leq \frac{1}{2}\Bigg[ |B^\prime (z)| - \frac{n|R(z)|^2}{\big\{M(R,1)\big\}^2} + \frac{2|R(z)|^2}{\big\{M(R,1)\big\}^2} \sum_{j=1}^{m}\frac{1}{1+|b_j|} \Bigg]M(R,1)\\
	&= \frac{1}{2}\Bigg[ |B^\prime (z)| - \frac{n|R(z)|^2}{\big\{M(R,1)\big\}^2}\bigg( 1 - \frac{2}{k+1} \bigg) -\frac{n|R(z)|^2}{\big\{M(R,1)\big\}^2}\frac{2}{k+1} + \frac{2|R(z)|^2}{\big\{M(R,1)\big\}^2} \sum_{j=1}^{m}\frac{1}{1+|b_j|} \Bigg]M(R,1)\\
	&=\frac{1}{2}\Bigg[ |B^\prime (z)| - \frac{n|R(z)|^2(k-1)}{(k+1)\big\{M(R,1)\big\}^2}  - \frac{2|R(z)|^2}{\big\{M(R,1)\big\}^2}\bigg\{\frac{n}{k+1}-\sum_{j=1}^{m}\frac{1}{1+|b_j|} \bigg\}\Bigg]M(R,1).
	\end{align*}
	This completes the proof of Theorem \ref{tR2}
\end{proof}

\begin{proof}[\textbf{Proof of Theorem \ref{tR1}}]
Inequality \eqref{tR1e3} together with Lemma \ref{lR4}, yields 
\begin{align*}
2\big| R^\prime (z) \big|^2 + n|B^\prime (z)||R(z)|^2 -2|B^\prime (z)| |R(z)|^2 \sum_{j=1}^{m}\frac{1}{1+|b_j|}&\leq \big|R^\prime (z)\big|^2 + \big|\big(R^* (z)\big)^\prime\big|^2\\
\nonumber &\leq \frac{1}{2}|B^\prime (z)|^2\big( M_1 ^2 + M_2 ^2 \big).
\end{align*}
Which can also be written as
\begin{align*}
4\big| R^\prime (z) \big|^2 
&\leq |B^\prime (z)|^2\big( M_1 ^2 + M_2 ^2 \big) + 4|B^\prime (z)| |R(z)|^2\sum_{j=1}^{m}\frac{1}{1+|b_j|}  - 2n|B^\prime (z)||R(z)|^2 \\
&= \Bigg[ |B^\prime (z)|^2 -\frac{2n(k-1)}{k+1} \frac{|R(z)|^2|B^\prime (z)|}{\big( M_1 ^2 + M_2 ^2 \big)} -  \frac{4|R(z)|^2|B^\prime (z)|}{\big( M_1 ^2 + M_2 ^2 \big)}\bigg\{ \frac{n}{1+k} - \sum_{j=1}^{m}\frac{1}{1+|b_j|} \bigg\} \Bigg]\big( M_1 ^2 + M_2 ^2 \big).
\end{align*}
Equivalently
\begin{align*}
\big| R^\prime (z) \big| \leq \frac{1}{2}\Bigg[ |B^\prime (z)|^2 -\frac{2n(k-1)}{k+1} \frac{|R(z)|^2|B^\prime (z)|}{\big( M_1 ^2 + M_2 ^2 \big)} -  \frac{4|R(z)|^2|B^\prime (z)|}{\big( M_1 ^2 + M_2 ^2 \big)}\bigg\{ \frac{n}{1+k} - \sum_{j=1}^{m}\frac{1}{1+|b_j|} \bigg\} \Bigg]^{\frac{1}{2}}\big( M_1 ^2 + M_2 ^2 \big)^{\frac{1}{2}}.
\end{align*}
Which being the desired conclusion, completes the proof of Theorem \ref{tR1}  
	\end{proof}
\begin{proof}[\textbf{Proof of Theorem \ref{tR3}}]
In the light of Lemma \ref{lR5}, one can easily see for $z\in\mathbb{J}$ and $R(z)$=$\frac{P(z)}{W(z)}$ with $P(z)=\alpha_m\prod_{j=1}^{m}(z-b_j)$, $\alpha_m \neq 0,$ $m\leq n, |b_j| \leq k\leq 1,$ for $ j=1,2,\dots,m$, that 
\begin{align*}
\nonumber \Re\bigg( \frac{zR^\prime (z)}{R(z)}\bigg) &= \Re\bigg( \frac{zP^\prime (z)}{P(z)} - \frac{zW^\prime (z)}{W(z)}  \bigg)\\
\nonumber &= \Re\bigg( \frac{zP^\prime (z)}{P(z)}  \bigg) - \Re\bigg( \frac{zW^\prime (z)}{W(z)}  \bigg)\\
\nonumber &= \sum_{j=1}^{m}\Re\bigg( \frac{z}{z-b_j} \bigg) - \frac{n-|B^\prime (z)|}{2}\\
& \geq \sum_{j=1}^{m}\frac{1}{1+|b_j|} - \frac{n-|B^\prime (z)|}{2}\\
&= \frac{|B^\prime (z)|}{2} - \frac{n}{2} + \sum_{j=1}^{m}\frac{1}{1+|b_j|}\\
&= \frac{|B^\prime (z)|}{2} + \frac{2m-n(1+k)}{2(1+k)} + \sum_{j=1}^{m}\frac{1}{1+|b_j|}-\frac{m}{1+k}.
\end{align*}
This for $z\in \mathbb{J}$, straightforwardly gives
\begin{align*}
\bigg| \frac{R^\prime (z)}{R(z)} \bigg|= \bigg| \frac{zR^\prime (z)}{R(z)} \bigg| \geq  \Re\bigg( \frac{zR^\prime (z)}{R(z)}\bigg) \geq \frac{|B^\prime (z)|}{2} + \frac{2m-n(1+k)}{2(1+k)} + \sum_{j=1}^{m}\frac{1}{1+|b_j|}-\frac{m}{1+k}.
\end{align*}
Which being equivalent to the desired inequality, completes the proof of Theorem \ref{tR3}
\end{proof}

\end{document}